\documentclass[11pt]{amsart}

\usepackage{amssymb}

\input xy

\xyoption{all}

\theoremstyle{plain}

\newtheorem{theorem}{Theorem}[section]

\newtheorem{corollary}[theorem]{Corollary}

\newtheorem{lemma}[theorem]{Lemma}

\newtheorem{conjecture}[theorem]{Conjecture}

\newtheorem{proposition}[theorem]{Proposition}

\newtheorem{fact}[theorem]{Fact}

\newtheorem{claim}[theorem]{Claim}

\theoremstyle{definition}

\newtheorem{definition}[theorem]{Definition}

\newtheorem{remark}[theorem]{Remark}

\newtheorem{notation}[theorem]{Notation}

\theoremstyle{remark}

\newcommand{\Aut}{\operatorname{Aut}}

\newcommand{\tp}{\operatorname{ga-tp}}

\newcommand{\gaS}{\operatorname{ga-S}}

\newcommand{\Hanf}{\operatorname{Hanf}}

\newcommand{\conc}{\Hat{\ }}

\newcommand{\sq}[2]{\sideset{^{#1}}{}{\operatorname{#2}}}

\newcommand{\lan}{\operatorname{L}}

\renewcommand{\phi}{\varphi}

\newcommand{\Union}{\bigcup}

\newcommand{\initial}\lessdot

\newcommand{\K}{\operatorname{\mathcal{K}}}

\newcommand{\id}{\operatorname{id}}

\newcommand{\LS}{\operatorname{LS}}

\newcommand{\C}{\mathfrak C}

\def\?{?\vadjust

{\vbox to 0pt{\vskip-7pt\hbox to 1.1\hsize{\hfill\huge ?!}}}}

\begin{document}

\title[Categoricity in Tame AECS]{Categoricity from one successor cardinal
in \\Tame Abstract Elementary Classes}

\author{Rami Grossberg}
\email[Rami Grossberg]{rami@andrew.cmu.edu}
\address{Department of Mathematics\\
Carnegie Mellon University\\
Pittsburgh PA 15213}

\author{Monica VanDieren}
\email[Monica VanDieren]{mvd@umich.edu}
\address{Department of Mathematics\\
University of Michigan\\
Ann Arbor MI 48109-1109}


\date{September 30, 2005}

\maketitle

\begin{abstract}
We prove that from categoricity in $\lambda^+$ we can get
categoricity in all cardinals $\geq\lambda^+$ in a $\chi$-tame
abstract elementary classes which has arbitrarily large models and
satisfies the amalgamation and joint embedding properties, provided
$\lambda>\LS(\K)$ and $\lambda\geq\chi$.

For the missing case when $\lambda=\LS(\K)$, we prove
that $\K$ is totally categorical provided that $\K$ is categorical in
$\LS(\K)$ and $\LS(\K)^+$.
\end{abstract}

\section{introduction}

The benchmark of progress in the development of a model theory for
abstract elementary classes (AECs) is Shelah's Categoricity Conjecture.

\begin{conjecture}\label{shelah's cat conj}
Let
$\K$ be an abstract elementary class.  If
$\K$ is categorical in some $\lambda >
\Hanf(\K)\footnote{$\Hanf(\K)$ is bounded
above by $\beth_{(2^{\LS(\K)})^+}$.  For a  definition of
$\Hanf(\K)$ see \cite{Gr1}.}$,
 then for every
$\mu\geq \Hanf(\K)$, $\K$ is categorical in $\mu$.
\end{conjecture}

With the exception of \cite{MaSh}, \cite{KoSh}, \cite{Sh 576}, \cite{ShVi}
and
\cite{Va} in which extra set theoretic assumptions are made, all work
towards Shelah's Categoricity Conjecture has taken place under the
assumption of the amalgamation property.  An AEC satisfies the
\emph{amalgamation property} if for every triple of models $M_0,M_1,M_2$
in which $M_0\prec_{\K}M_1$ and $M_0\prec_{\K}M_2$ there exist
$\K$-mappings $g_1$ and $g_2$ and an amalgam $N\in\K$ such that the
diagram below commutes.

\[
\xymatrix{\ar @{} [dr] M_1
\ar[r]^{g_1}  &N \\
M_0 \ar[u]^{\id} \ar[r]_{\id}
& M_2 \ar[u]_{g_2}
}
\]

Under the assumption of the amalgamation property, there is a natural
generalization of first order types.  However, types are no longer
identified by consistent sets of formulas.
Since we assume the amalgamation and joint embedding
properties, we may work inside a large monster model which we denote by
$\C$.  We use the notation $\Aut_M(\C)$ to represent the group of
automorphisms of $\C$ which fix $M$ pointwise.  With the
amalgamation property, we can define the Galois-type of an element $a$
over a model
$M$, written
$\tp(a/M)$.  We say two elements
$a,b\in\C$ realize the same Galois-type over a model $M$ iff there is an
automorphism $f$ of $\C$ such that $f(a)=b$ and $f\restriction M=\id_M$.
We abbreviate the set of all Galois-types over a model $M$ by $\gaS(M)$.
An AEC is
\emph{Galois-stable} is $\mu$ if for every model $M$ of $\K$ of
cardinality
$\mu$, there are only $\mu$ many Galois-types over $M$.  See \cite{Gr1} or
\cite{Ba1} for a survey of the development of these concepts.

In the first author's Ph.D. thesis and \cite{GrVa2}, we isolated the
notion of tameness in order to develop a stability theory for a wide
spectrum of non-elementary classes.  An abstract elementary class
satisfying the amalgamation property is said to be \emph{$\chi$-tame} if
for every model
$M$ in $\K$ of cardinality
$\geq\chi$ and every $p\neq q\in\gaS(M)$, there is a submodel $N$ of $M$
of cardinality $\chi$ such that $p\restriction N\neq q\restriction N$.
A class $\K$ is said to be \emph{tame} if it is $\chi$-tame for some
$\chi$. In other words, tameness captures the local character of
consistency.

\emph{All} families of AECs that are known to  have a
structural theory satisfy the amalgamation property and  are tame
\footnote{While there are structural results for continuous model theory, this
context is
not an AEC.  The classification theory for continuous model theory is
parallel to the Buechler-Lessmann paper on homogeneous models
\cite{BuLe}.  One could apply the definition of tame to classes satisfying
the same properties as models of a continuous theory.  In this view,
continuous model theory is tame.}.
In fact several examples of tame class fail to be homogeneous or even
excellent.
\begin{enumerate}
\item Elementary classes.

\item
Homogeneous model theory (as Galois-types are sets of formulas).

\item
The class of atomic models of a first-order theory (from
\cite{Sh 87a}). I.e. the class  introduced to study the spectrum
function of $L_{\omega_1,\omega}$ sentence (under mild assumptions) is an
example of a tame AEC.

\item
Let $\K$ be an AEC,   and suppose there exists \emph{  $\kappa$  strongly
compact cardinal}
such that $\LS(\mathcal{K})<\kappa$.  Let $\mu_0:=\beth_{(2^\kappa)^+}$.
Makkai and Shelah prove that if
$\K$ is categorical in some
      $\lambda^{+}>\mu_0$ then has the
AP.  By further results of
      \cite{MaSh} the Galois-types can be identified with sets of formulas
from $L_{\kappa,\kappa}$.  Thus
      $\K$ is $\kappa$-tame.

\item
The class of algebraically closed fields with pseudo-exponentiation
studied by Zilber is tame.

\item  Using the method of \cite{GrKo} Villaveces and Zambrano in \cite{ViZa}
have shown that
the class of Hrushovski's fusion $\K_{fus}$
is $\aleph_0$-tame.

\item Baldwin \cite{Ba2} combining arguments from  \cite{GrKo} and
\cite{Zi2} have shown that the class $\K$ of two sorted structures
$(V,A)$ when $A$ is  semi-abelian  with a group homomorphism $\exp$ from
a finite dimensional $\mathbf Q$-vector space $V$ onto $A$ with kernel
$\mathbf Z^N$ is $\aleph_0$-tame AEC with AP (details are in section 4 of
\cite{Ba2}).

\item It is a corollary of \cite{GrKo} that good frames that are excellent (in
the sense of \cite{Sh 705}) are tame.
\end{enumerate}

As further evidence to the importance of tame AECs,
recent
progress on Shelah's Categoricity Conjecture has been made under the
assumption of tameness by combining the work of \cite{Sh 394} with
\cite{GrVa1}.

\begin{fact}\label{sh and grva}
Suppose $\K$ is a $\chi$-tame abstract elementary class satisfying the
amalgamation and joint embedding  properties.
  Let
$\mu_0:=\Hanf(\K)$.  If $\chi\leq\beth_{(2^{\mu_0})^+}$ and
$\K$ is categorical in some
$\lambda^+>\beth_{(2^{\mu_0})^+}$, then $\K$ is categorical in $\mu$ for
all
$\mu>\beth_{(2^{\mu_0})^+}$.
\end{fact}

Previous results (e.g. \cite{Sh 87a}, \cite{Sh 87b}, \cite{MaSh},
\cite{KoSh}, \cite{Sh 472} and \cite{Sh 705}) of Shelah in the
direction of upward categoricity required not only model-theoretic
assumptions but also set-theoretic assumptions.  An interesting feature
of our work is that it is an upward categoricity transfer theorem in ZFC.
In particular it can be viewed as an improvement of the main result of
\cite{MaSh} where the assumption of existence of a strongly compact cardinal is made.

One distinction between
Fact
\ref{sh and grva} and Conjecture
\ref{shelah's cat conj} is that Fact \ref{sh and grva}  applies only to
classes which are categorical above the second Hanf number,
$\beth_{(2^{\Hanf(\K)})^+}$.  One motivation for this paper is to improve
Fact
\ref{sh and grva} getting a better approximation to Conjecture
\ref{shelah's cat conj} for tame abstract elementary classes.  In fact
our results extend beyond the scope of Conjecture \ref{shelah's cat
conj} since we are able, for instance,
to conclude that for a $\LS(\K)$-tame abstract elementary class with
arbitrarily large models satisfying the amalgamation and joint embedding
properties if the class is categorical in $\LS(\K)$ and $\LS(K)^+$ then
the class is categorical in all $\mu\geq\LS(\K)$.

In his paper \cite{Sh 394}, Shelah proved that from categoricity in
$\lambda^+$  above the second Hanf number, one could deduce categoricity
below
$\lambda^+$.  Under the additional assumption of tameness, we provide an
argument to transfer categoricity in
$\lambda^+$ upwards in
\cite{GrVa1}.  The main step in our proof is:

\begin{fact}[Corollary 4.3 of \cite{GrVa1}]\label{going up grva}
Suppose that $\K$ has arbitrarily large models, satisfies the
amalgamation property and is $\chi$-tame with
$\chi\geq\LS(\K)$. If
$\K$  is categorical in both $\lambda^+$ and $\lambda$
with $\lambda\geq\chi$ and $\lambda>\LS(\K)$, then
$\K$ is categorical in every
$\mu$ with $\mu\ge\lambda$.
\end{fact}

A  breakthrough  in \cite{GrVa1} was to go from categoricity
in
$\lambda^+$ to categoricity in $\lambda^{++}$ when $\lambda^+$ was above
the second $\Hanf$ number of the class.  Working under the assumption of
categoricity above the second $\Hanf$ number provided us the convenience
of categoricity in $\lambda$ with an application of
\cite{Sh 394}.

Recently, Lessmann expressed interest in whether or not the upward
categoricity transfer theorem (Fact \ref{going up grva}) could be proved
from categoricity in only one successor cardinal.  He communicated to us
that he could use our methods along with quasi-minimal types and countable
superlimits to prove the desired result for $\aleph_0$-tame classes with
$\LS(\K)=\aleph_0$
\cite{Le}, but was unable to prove it when $\LS(\K)$ is uncountable.  This paper answers
Lessmann's question. Using the ideas and arguments from
\cite{GrVa1} along with quasi-minimal types, we deduce from categoricity
in
$\lambda^+$   categoricity in $\lambda^{++}$ for $\lambda>\LS(\K)$ with no
restrictions on the size of $\LS(\K)$ or the tameness cardinal.  We also
improve Fact
\ref{going up grva} by removing the assumption that $\lambda>\LS(\K)$.

Our proof that categoricity in $\lambda^+$ implies categoricity in
$\lambda^{++}$ under the described setting involves showing that there
are nice minimal types (which we have called  deep-rooted quasi-minimal)
over limit models, and these quasi-minimal types have  no Vaughtian pairs
of cardinality
$\lambda^{++}$.  Then using a characterization of limit models (Theorem
4.1 from
\cite{GrVa1}), we show that this is enough to prove the model of
cardinality
$\lambda^{++}$ is saturated.

We are grateful to John Baldwin and Olivier Lessmann for
asking questions without which this paper would not exist.

\section{Preliminaries}
Throughout this paper, we make the assumptions that our abstract
elementary class
$\K$ has arbitrarily large models and satisfies the joint embedding and
amalgamation properties.    We will also assume that the class is
$\chi$-tame.  We let $\K_\mu$ stand for the set of all models
of $\K$ of cardinality $\mu$.  In the natural way, we use $\K_{\leq\mu}$
and $\K_{\geq\mu}$.  We will be using notation and definitions consistent
with \cite{GrVa1}. Many of the propositions can be proved in more general
settings, but we leave an exploration of those possibilities for future
work.


In abstract elementary classes saturated models
have various guises.  In some cases, it is more prudent to work with
a limit model as opposed to a saturated model.
\begin{definition}
\begin{enumerate}
\item \index{universal over!$\kappa$-universal over}

We say $N$ is \emph{universal over $M$} iff
for every $M'\in\K_{\|M\|}$ with $M\prec_{\K}M'$ there exists
a $\K$-embedding $g:M'\rightarrow N$ such that
$g\restriction M=\id_{M}$:

\[
\xymatrix{\ar @{} [dr] M'
\ar[drr]^{g}  & \\
M \ar[u]^{\id} \ar[rr]_{\id}
&& N
}
\]

\item For $M\in\K_\mu$,  $\sigma$ a limit ordinal with
$\sigma\leq \mu$ and $M'\in \K_{\mu}$
we say that $M'$ is a \emph{$(\mu,\sigma)$-limit over $M$}\index{limit
model!$(\mu,\sigma)$-limit model over $M$} iff there exists a
$\prec_{\K}$-increasing and continuous sequence of models $\langle M_i\in
\K_{\mu}\mid i<\sigma
\rangle$
such that
\begin{enumerate}

\item $M= M_0$,

\item $M'=\Union_{i<\sigma}M_i$
and

\item\label{univ cond in defn} $M_{i+1}$ is universal over $M_i$.

\end{enumerate}
\end{enumerate}
\end{definition}

While using back and forth one can show that any two
$(\mu,\sigma)$-limit models  are isomorphic
to show that all  $(\mu,\sigma_1)$-limit models are
isomorphic to a $(\mu,\sigma_2)$-limit model is not so obvious.
 We will be
using the following fact which is a consequence of
\cite{Va}; or see
\cite{GrVa3} for an exposition and proof.

\begin{fact}[Uniqueness of Limit Models]\label{unique limits}  Suppose
that $\K$ is an abstract elementary class satisfying the amalgamation
property and is categorical in some $\lambda$.  If $\LS(\K)<\mu<\lambda$
and
$M_0\in\K_\mu$, then for every two limit models
$M_1$ and $M_2$  over $M_0$, if $M_1$ and $M_2$ both have the same
cardinality $\kappa<\lambda$, then they are isomorphic over $M_0$.
\end{fact}

\begin{notation}
In light of Fact \ref{unique limits}, when the cardinality of the limit
model is clear,  we omit the parameters
$\mu$ and
$\sigma$ and refer to $(\mu,\sigma)$-limit models as \emph{limit models}.
\end{notation}

A corollary of Fact \ref{unique limits} is that
\begin{proposition}\label{limits are sat}
 Assuming categoricity in $\lambda$ and the joint embedding
and amalgamation properties, for $\mu$ with $\LS(\K)<\mu<\lambda$, every
saturated model of cardinality
$\mu$  is also a
$(\mu,\sigma)$-limit model for any limit ordinal $\sigma<\mu^+$.
\end{proposition}
\begin{proof}
First we show that any limit model of cardinality $\mu$ is saturated.
Then by our assumptions and the uniqueness of saturated
models (Lemma 0.26 of
\cite{Sh 576}), we can conclude that any saturated model of cardinality
$\mu$ is isomorphic to a limit model of cardinality $\mu$.

Suppose that $M$ is a limit model of cardinality $\mu$.  Fix $\kappa$
such that $\LS(\K)\leq\kappa<\mu$. Fix $N_1\prec_{\K}M$  of
cardinality
$\kappa$ and $p=\tp(a/N_1,N_2)$ with $\|N_2\|=\kappa$.
 Since $M$ is a limit model, we can
find a continuous decomposition of $M$ into $\langle M_i\mid
i<\kappa^+\rangle$ such that each $M_i$ is a model of cardinality $\mu$
and $M_{i+1}$ is universal over $M_i$.  By the regularity of
$\kappa^+$, we can find $i<\kappa^+$ such that $N_1\prec_{\K}M_i$.
 Invoking the amalgamation property, we can amalgamate
$N_2$ and
$M_i$ over $N_1$ as in the diagram below:

\[
\xymatrix{\ar @{} [dr] N_2
\ar[r]^{g}  &N^* \\
N_1 \ar[u]^{\id} \ar[r]_{\id}
& M_i \ar[u]^{\id}
}
\]

We may assume that the amalgam $N^*$ has cardinality $\mu$.  Since
$M_{i+1}$ is universal over $M_i$ we can extend the commutative diagram:

\[
\xymatrix{\ar @{} [dr] N_2
\ar[r]^{g}  &N^* \ar[dr]^{f}& \\
N_1 \ar[u]^{\id} \ar[r]_{\id}
& M_i \ar[u]^{\id} \ar[r]_{\id}&M_{i+1}
}
\]

Notice that $f\circ g$ witnesses that $f(g(a))\in M_{i+1}$ realizes
$\tp(a/N_1,N_2)$.

\end{proof}

Galois-stability and the
amalgamation property are enough to establish the existence of limit
models (see
\cite{Sh 600} for the statement and \cite{GrVa2} for a
proof).  Limit models exist in categorical AECs since
categoricity implies Galois-stability:
\begin{fact}[Claim 1.7(a) of \cite{Sh 394} or see
\cite{Ba1} for a proof]\label{cat implies stable}

If $\K$ is categorical in $\lambda> \LS(\K)$, then $\K$ is
Galois-stable in all $\mu$ with $\LS(\K)\leq\mu<\lambda$.
\end{fact}

Another consequence of $\mu$-stability is the existence of minimal types.
As a replacement for first order strongly minimal types, Shelah has
suggested using minimal types in \cite{Sh 394}.  We in \cite{GrVa1}
found that a more restrictive minimality condition (rooted minimal) could
be used to transfer categoricity upward.
\begin{definition}\label{strongly rooted min defn}
Let $M\in\K$ and $p\in\gaS(M)$ be given.
\begin{enumerate}
\item $p$ is said to be \emph{minimal} if it is both non-algebraic (that
is, it is not realized in $M$) and for any $N\in\K$ extending $M$ there is
at most one non-algebraic extension of $p$ to $N$.
\item A minimal type $p$ is said to be \emph{rooted minimal} iff there is
some $M_0\prec_{\K}M$ with $M_0\in\K_{<\|M\|}$ such that $p\restriction
M_0$ is also minimal.  $M_0$ is called a \emph{root} of $p$.

\end{enumerate}
\end{definition}

\begin{fact}[Density of Minimal Types \cite{Sh 394}]\label{density exist
minimal} Let $\mu>\LS(\K)$.
If $\K$ is Galois-stable in $\mu$, then for every $N\in\K_\mu$
and every
$q\in\gaS(N)$, there are $M\in\K_\mu$ and $p\in\gaS(M)$ such that
$N\preceq_{\K}M$, $q\leq p$ and $p$ is minimal.
\end{fact}

The main obstacle of minimal types in this context is that while there
are minimal types in stable AECs, the minimal types may be trivially
minimal, meaning that the minimal type has \emph{no} non-algebraic
extensions.  As in \cite{Sh 48} and \cite{Zi} we replace this notion of
minimality with quasi-minimality.

Since a non-algebraic
type may not have any non-algebraic extensions, we distinguish these
non-algebraic types from the well-behaved non-algebraic types.  A
non-algebraic type
$p\in\gaS(M)$ is
\emph{big} iff for every $M'\succ_{\K} M$ of cardinality $\|M\|$, there is a
non-algebraic extensions of $p$ to $M'$ (see Definition 6.1 of \cite{Sh
48}).  Notice that this is equivalent to requiring that there is a big
extension of $p$ to
$M'$.

Almost thirty years after Shelah's \cite{Sh 48}, Zilber rediscovered the
notion of minimality and used perhaps the better notation
quasi-minimality to distinguish it from the first order relatives. As in
\cite{Zi}, we say a big type
$p$ is
\emph{quasi-minimal} iff for  any $N\in\K$ extending $M$ there is
at most one non-algebraic extension of $p$ to $N$.    Analogous to the
minimal case, we can define
\emph{deep-rooted quasi-minimal}. Most of the results concerning minimal
types can be proved for quasi-minimal types with minimal work.

We will show that quasi-minimal types exist in Section \ref{s:exist
quasiminimal}.  For now notice that the assumptions of the amalgamation
property and no maximal models give us the following:

\begin{remark}\label{big exist}
For any $M\in\K$, there exists
$p\in\gaS(M)$ such that
$p$ is big.
\end{remark}

Another consequence of $\mu$-stability is that $\mu$-splitting is
well-behaved and the notions of non-algebraic and big types over limit
models are the same.  We begin by reviewing some basic facts about
$\mu$-splitting.

 For $M\in\K_{\geq\mu}$ and $N\prec_{\K}M$ we
say that
$p\in\gaS(M)$
\emph{$\mu$-splits over $N$} iff there exist
two models $N_1,N_2\in\K_\mu$ and an isomorphism $f:N_1\cong N_2$ such
that
$N\prec_{\K}N_l\prec_{\K}M$ for $l=1,2$; $f\restriction N=\id_N$ and
$p\restriction N_2\neq f(p\restriction N_1)$.  Under the assumption of
categoricity, $\mu$-splitting has an extension property (See Corollary 2
of \cite{Ba2} or Theorem 12.8 of \cite{Ba1}) in addition to the existence
property which follows from Galois-stability in
$\mu$ (see Lemma 6.3 of \cite{Sh 394}):

\begin{fact}\label{Baldwin
non-split ext} Suppose that $\K$ is categorical
in some
$\lambda>\LS(\K)$.
Let $\mu$ be a
cardinal such that
$\LS(\K)\leq\mu$ and let $\sigma$ be a limit ordinal with
$\LS(\K)\leq\sigma<\mu^+$.
Then, for every
$(\mu,\sigma)$-limit model $M$ and every type $p\in \gaS(M)$, there exists
 $N\precneqq_{\K}M$ of cardinality $\mu$ such that for
every
$M'\in\K_{\leq\lambda}$ extending $M$, there exists $q\in \gaS(M')$ an
extension of $p$ such
that $q$ does not $\mu$-split over $N$.  In particular $p$ does not
$\mu$-split over $N$.

Moreover, if $M$ is a $(\mu,\sigma)$-limit model witnessed by $\langle
M_i\mid i<\sigma\rangle$, then there is a $i<\sigma$ such that $p$ does
not $\mu$-split over $M_i$.
\end{fact}

The only other property of $\mu$-splitting that we will explicitly use is
an observation that non-splitting extensions of non-algebraic types
remain non-algebraic.
\begin{fact}[Corollary 2.8 of \cite{GrVa1}]\label{non
split ext of non alg}
Let $N,M,M'\in\K_\mu$ be
such that
$M'$ is universal over
$M$ and
$M$ is a limit model over $N$.
Suppose that $p\in\gaS(M)$ does not $\mu$-split over $N$ and $p$ is
non-algebraic.  For every $M'\in\K$ extending $M$ of cardinality
$\mu$, if $q\in \gaS(M')$ is an extension of $p$ and does not $\mu$-split
over $N$,  then $q$ is non-algebraic.
\end{fact}

We can use non-splitting to show that
\begin{fact}\label{big is non-alg}
Suppose that $\K$ is categorical
in some
$\lambda>\LS(\K)$ and $\mu$ is a cardinal $<\lambda$.
 If $M$ is a limit model of cardinality $\mu$,
then $p\in\gaS(M)$ is non-algebraic iff $p$ is big.
\end{fact}

\begin{proof}As in the proof of Theorem I.4.10 \cite{Va} or see
Proposition 1.16\cite{Le}.  At the referee's suggestion, we have included
a proof here.

Clearly every big type is non-algebraic.  Suppose $M$ is a limit model
witnessed by $\langle M_i\mid i<\sigma\rangle$ and $p=\tp(a/M)$
is non-algebraic. By Fact \ref{Baldwin
non-split ext}, there is an
$i<\sigma$ such that
$p$ does not $\mu$-split over $M_i$.

Let $M'$ be a $\K$-extension of $M$ of cardinality $\mu$.  We now show
that
$p$ can be extended to a non-algebraic type $p'\in\gaS(M')$.  By the
definition of limit  model and our choice of $\langle M_i\mid
i<\sigma\rangle$, we know that $M_{i+1}$ is universal over $M_i$.  Thus
there is a $\K$-mapping $h':M'\rightarrow M_{i+1}$ such
that
$h\restriction M_i=\id_{M_i}$.  Because of we are working inside a
monster model, we can extend $h'$ to $h\in\Aut_{M_i}(\C)$.
Our candidate for a non-algebraic extension of $p$ to $M'$ will be
$p':=\tp(h^{-1}(a)/M')$.  Immediately we see that $p'$ is non-algebraic
since $\tp(a/h(M'))$ was non-algebraic.

 We claim that $p'$ is in fact an extension of $p$, that is that
$\tp(h^{-1}(a)/M)=\tp(a/M)$.
 By
monotonicity, of non-splitting, we have that
$\tp(a/h(M')$ does not $\mu$-split over $M_i$.
By invariance, we have $\tp(h^{-1}(a)/M')$ also does not $\mu$-split
over $M_i$.  Now if $\tp(h^{-1}(a)/M)\neq \tp(a/M)$, we would
witness that
$\tp(h^{-1}(a)/M')$ $\mu$-splits over $M_i$ via the mapping $h$.  Thus
$p'\restriction M=p$ as required.

\end{proof}

We now go into some details of a common construction in AECs.  A
variation of the proposition appears in the literature as Claim 0.31(2)
of \cite{Sh 576} and in the proof of Theorem II.7.1 of
\cite{Va}, we isolate it here as Lemma \ref{Baldwin's coherence}.
After detailing Lemma \ref{Baldwin's coherence} to John Baldwin in e-mail
correspondence in the Summer of 2004, we decided to include the proof
here.

\begin{lemma}\label{Baldwin's coherence}
Suppose $\langle M_i\mid i<\alpha\rangle$ is an $\prec_{\K}$-increasing
and continuous chain of models. Further assume that $\langle
p_i\in\gaS(M_i)\mid i<\alpha\rangle$ is an increasing chain of types such
that there are $a_i\in\C$ with
$a_i\models p_i$ and
$\prec_{\K}$-mappings $f_{i,j}\in\Aut_{M_i}(\C)$ with $f_{i,j}(a_i)=a_j$
for $i\leq j<\alpha$ such that for $i\leq j\leq k$ we have that
$f_{i,k}=f_{j,k}\circ f_{i,j}$.  Then there exists $a_\alpha\in\C$
realizing each $p_i$ and there are $f_{i,\alpha}\in\Aut_{M_i}(\C)$ with
$f_{i,\alpha}(a_i)=a_\alpha$.
\end{lemma}

The proof uses direct limits, so we will review some facts first.
Using the axioms of AEC and Shelah's Presentation Theorem, one can show
that the union axiom of the definition of AEC has an alternative
formulation (see
\cite{Sh 88} or Chapter 16 of
\cite{Gr2}):

\begin{definition}\index{directed set}A partially ordered set $(I,\leq)$
is
\emph{directed} iff for every $a,b\in I$, there exists
$c\in I$ such that $a\leq c$ and $b\leq c$.
\end{definition}

\begin{fact}[P.M. Cohn 1965]\label{direct limits}\index{direct
limit!existence of in AECs} Let $(I,\leq)$ be a directed set.  If $\langle
M_t\mid t\in I\rangle$ and $\{h_{t,r}\mid t\leq r\in I\}$ are such that
\begin{enumerate}
\item for $t\in I$, $M_t\in \K$
\item for $t\leq r\in I$, $h_{t,r}:M_t\rightarrow M_r$ is a
$\prec_{\K}$-embedding and
\item for $t_1\leq t_2\leq t_3\in I$,
$h_{t_1,t_3}=h_{t_2,t_3}\circ h_{t_1,t_2}$ and $h_{t,t}=\id_{M_t}$,
\end{enumerate}
then, whenever $s=\lim_{t\in I} t$, there exist
$M_s\in \K$ and $\prec_{\K}$-mappings
$\{h_{t,s}\mid t\in I\}$ such that
$$h_{t,s}:M_t\rightarrow M_s, M_s=\Union_{t<s}h_{t,s}(M_t)\text{ and}$$
$$\text{for }t_1\leq t_2\leq s, h_{t_1,s}=h_{t_2,s}\circ h_{t_1,t_2}
\text{ and }h_{s,s}=id_{M_s}.$$
\end{fact}

\begin{remark}
Cohn's proof gives us that $M_s$ is an $\lan(\K)$-structure.  To show that
$M_s\in \K$ and that $h_{t,s}$ are $\K$-embeddings we use Shelah's
presentation theorem.
\end{remark}

\begin{definition}
\begin{enumerate}

\item\index{directed system} $(\langle M_t\mid t\in I\rangle,\{h_{t,s}\mid
t\leq s\in I\})$ from Fact \ref{direct limits} is called a \emph{directed
system}.
\item\index{direct limit!definition}
We say that $M_s$ together with $\langle h_{t,s}\mid t\leq s\rangle$
satisfying the conclusion of Fact \ref{direct limits} is \emph{a
direct limit of $(\langle M_t\mid t<s\rangle,\{h_{t,r}\mid t\leq r<s\})$}.
\end{enumerate}
\end{definition}

\begin{proof}[Proof of Lemma \ref{Baldwin's coherence}]
Let $\langle p_i\in\gaS(M_i)\mid i<\alpha\rangle$ be an increasing chain
of types and
$\langle M_i\mid i<\alpha\rangle$ a $\prec_{\K}$-increasing chain of
models with $\langle f_{i,j}\mid i\leq j<\alpha\rangle$ and $a_i$ as in
the statement of the lemma.
Notice that $(\langle \C_i\mid
i<\alpha\rangle,\langle f_{i,j}\mid i\leq j<\alpha\rangle)$ forms a
directed system where $\C_i=\C$ for all $i$.  Let $\C^*_\alpha$ and
$\langle f^*_{i,\alpha}\mid i\leq\alpha\rangle$ be a direct limit to this
system.   Outright we don't have much control over this limit, but by the
following claims we will be able to chose a limit $(\C_\alpha,\langle
f_{i,\alpha}\mid i\leq\alpha\rangle)$ so that
$\Union_{i<\alpha}M_i\preceq_{\K}\C_\alpha=\C$ and
$f_{i,\alpha}\restriction M_i=\id_{M_i}$.

First notice that we can take $\C_\alpha$ to be $\C$ since
 a direct limit of automorphisms is an
isomorphism using the construction of direct limits from \cite{Gra}.

\begin{claim}\label{inc claim}
$\langle f^*_{i,\alpha}\restriction M_i\mid i\leq\alpha\rangle$ is
increasing.
\end{claim}

\begin{proof}[Proof of Claim \ref{inc claim}]
Let $i< j<\alpha$ be given.  By construction
$$f_{i,j}\restriction M_i=\id_{M_i}.$$
An application of $f^{*}_{j,\alpha}$ yields
$$f^{*}_{j,\alpha}\circ f_{i,j}\restriction
M_i=f^{*}_{j,\alpha}\restriction M_i.$$
Since $f^{*}_{i,\alpha}$ and $f^{*}_{j,\alpha}$ come from a direct
limit of the system which includes the mapping $f_{i,j}$,
we have
$$
f^{*}_{i,\alpha}\restriction M_i=
f^{*}_{j,\alpha}\circ f_{i,j}\restriction
M_i.$$
Combining the equalities yields
$$  f^{*}_{i,\alpha}\restriction
M_i=f^{*}_{j,\alpha}\restriction M_i.$$
This completes the proof of Claim \ref{inc claim}.

\end{proof}

By the claim, we have that
$f:=\Union_{i<\alpha}f^{*}_{i,\alpha}\restriction M_i$ is a
$\prec_{\K}$-mapping from
$\Union_{i<\alpha}M_i$ onto
$\Union_{i<\alpha}f^{*}_{i,\alpha}(M_i)$.  Since $\C$ is saturated and
model homogeneous,  we can  extend $f$
to
$F\in \Aut(\C)$.

Now
consider the direct limit defined by
$\C_\alpha:=F^{-1}(\C^{*}_\alpha)$ with
$\langle f_{i,\alpha}:=F^{-1}\circ f^{*}_{i,\alpha}\mid i<
\alpha\rangle$ and $f_{\alpha,\alpha}=\id_{\C_\alpha}$.
 Notice that
$f_{i,\alpha}\restriction M_i=F^{-1}\circ
f^{*}_{i,\alpha}\restriction M_i=\id_{M_i}$ for
$i<\alpha$.  Thus $\Union_{i<\alpha}M_i\preceq_{\K}\C_\alpha$.

Let $a_\alpha:=f_{0,\alpha}(a_0)$.  The following argument explains why
$\tp(a_\alpha/\Union_{i<\alpha}M_i)$ is an upper bound for $\langle
p_i\mid i<\alpha\rangle$.

\begin{claim}\label{inc type claim}$\tp(a_\alpha/M_i)=\tp(a_i/M_i)$ for
all
$i<\alpha$.
\end{claim}
\begin{proof}[Proof of Claim \ref{inc type claim}]
Fix $i<\alpha$.  Notice that by the definition of direct limit we have
$a_\alpha=f_{0,\alpha}(a_0)=f_{i,\alpha}\circ f_{0,i}(a_0)$.  But by our
choice of $f_{0,i}$ we know that $f_{0,i}(a_0)$ is actually $a_i$.  Thus
$f_{i,\alpha}$ is an automorphism of $\C$ fixing $M_i$ taking $a_i$ to
$a_\alpha$.  So $\tp(a_i/M_i)$ and $\tp(a_\alpha/M_i)$ must be the same.

\end{proof}

\end{proof}

\begin{lemma}\label{can find the right auto}
Suppose $\langle M_i\mid i<\alpha\rangle$ is an $\prec_{\K}$-increasing
and continuous chain of limit models. If $\langle p_i\in\gaS(M_i)\mid
i<\alpha\rangle$ is an increasing chain of quasi-minimal types and
$\alpha$ is a limit ordinal, then  we can find $a_i\in\C$ with $a_i\models p_i$
and
$\prec_{\K}$-mappings $f_{i,j}\in\Aut_{M_i}(\C)$ with $f_{i,j}(a_i)=a_j$
for $i\leq j<\alpha$ such that for $i\leq j\leq k$ we have that
$f_{i,k}=f_{j,k}\circ f_{i,j}$.
\end{lemma}
\begin{proof}
We find $a_i$ and $f_{k,i}$ by induction on $i$.  For $i=0$, take
$a_0\in\C$ to be some realization of $p_0$ and $f_{0,0}:=\id_{\C}$.
Suppose that we have defined $a_i$ and $f_{k,i}$ for all $k\leq i$.  Let
$a_{i+1}$ be some realization of $p_{i+1}$ in $\C$.  Since the types are
increasing, we can find $f_{i,i+1}\in\Aut_{M_i}\C$ with
$f_{i,i+1}(a_i)=a_{i+1}$.  Define $f_{k,i+1}:=f_{i,i+1}\circ f_{i,i+1}$.
We use quasi-minimal types to get past  limit stage.  Suppose that we have
defined
$f_{j,k}$ for all $j\leq k<i$ with $i$  a limit ordinal.   By Lemma
\ref{Baldwin's coherence} there exists $a^*\in\C$ and
$f_{j,i}\in\Aut_{M_j}\C$ with $f_{j,i}\restriction M_j=id_{M_j}$ and
$f^*_{j,i}(a_j)=a^*$.  This $a^*$ comes from a direct limit construction
and may not realize the same type as $a_i$ over $M_i$.  However,
$\tp(a^*/M_i)$ is a non-algebraic extension of $\tp(a_0/M_0)$, which was
quasi-minimal.  Since $M_i$ is also a limit model, then $\tp(a^*/M_i)$
is big.  So, we can actually conclude, by quasi-minimality that  the types
of
$a^*$ and $a_i$ over $M_i$ agree.  So we can fix $g\in\Aut_{M_i}(\C)$
such that $g(a^*)=a_i$.  Then $f_{j,i}:=g\circ f^*_{j,i}$ is as required.
\end{proof}

\begin{corollary}\label{non-alg limits exist}
Suppose $\langle M_i\mid i<\alpha\rangle$ is an $\prec_{\K}$-increasing
and continuous chain of limit models. If $\langle p_i\in\gaS(M_i)\mid
i<\alpha\rangle$ is an increasing chain of quasi-minimal types and
$\alpha$ is a limit ordinal, then
 there is a
$p_\alpha\in\gaS(\Union_{i<\alpha}M_i)$ extending each of the $p_i$.

\end{corollary}

\begin{proof}
Follows from Lemma \ref{can find the
right auto} and Lemma
\ref{Baldwin's coherence}.

\end{proof}

\section{Deep-rooted minimal types}\label{s:exist quasiminimal}

 The main aim of this section is to prove the existence of deep-rooted
quasi-minimal types.  We will use the idea of Shelah's density of minimal
types to do this.  Our work generalizes Lemma 6.6 of \cite{Sh 48} where
Shelah proves the existence of (quasi)-minimal types using a rank
function.

First notice that if the class is tame, then any big
extension of a quasi-minimal type  is also quasi-minimal:

\begin{proposition}[Monotonicity of Minimal Types]\label{monotonicity of
minimal types} Suppose
$\K$ is
$\chi$-tame for some
$\chi$ with $\mu\geq\chi$.
If $p\in\gaS(M)$ is quasi-minimal with $M\in\K_\mu$, then for all
$N\in\K$ extending $M$ and every $q\in \gaS(N)$ extending $p$,
if
$q$ is big then $q$ is quasi-minimal.  If $N$ is a limit model, then the
assumption that $q$ is big can be replaced with non-algebraic.
\end{proposition}

\begin{proof}
The last sentence of the claim is Proposition 2.2 of \cite{GrVa1} once we
notice Fact \ref{big is non-alg}. The proof of the rest of the claim is
similar, but  we include details here for completeness.

Suppose that $p$ is a quasi-minimal type over $M$ with a big extension
$q$ to $N$.  For the sake of contradiction assume that $q$ is not
quasi-minimal.  Then there exist distinct $q_1$ and $q_2$ non-algebraic
extensions of $q$ to some model $N'$.  By tameness,
there exists $M'$ of cardinality $\mu$ such that $M\prec_{\K}M'$ and
$q_1\restriction M'
\neq q_2\restriction M'$.
Since $q_1\restriction M'$ and $ q_2\restriction M'$ are both
non-algebraic
extensions of $p$ we have a contradiction to the quasi-minimality of $p$.
\end{proof}


%

Similar to the proof of the density of minimal types, Fact \ref{density
exist minimal}, we get quasi-minimal types.  Moreover, instead of a
density result, we can actually find quasi-minimal types over every limit
model.  This is one of the obstacles in working with minimal types.

\begin{proposition} [Existence of Quasi-Minimal Types over
Limits]\label{exist minimal}Suppose $\K$ is Galois-stable in $\mu$ and
$M\in\K_\mu$ is a limit model.  Then there exists a quasi-minimal type
over $M$.
\end{proposition}

\begin{proof}
 We build a tree of types, but restrict ourselves to
limit models throughout the construction.
Suppose for the sake of contradiction that $M\in\K_\mu$ is a limit model
and that there are no quasi-minimal types over $M$.
By Remark \ref{big exist} we can
fix $p\in\gaS(M)$ a
big type.
 By induction on $i<\mu^+$ we
build a $\prec_{\K}$-increasing and continuous chain of models, $\langle
M_i\mid i<\mu^+\rangle$ and a tree of types $\langle
p_\eta\mid\eta\in\sq{<\mu^+}2\rangle$
satisfying
\begin{enumerate}
\item $M_i$ is a limit model of cardinality $\mu$
\item\label{get limit} $M_{i+1}$ is a limit model over $M_i$
\item $p_\eta=\tp(a_\eta/M_i)$ is big where $i$ is
the length of
$\eta$
\item $p_{\eta\conc\langle 0\rangle}\neq p_{\eta\conc\langle 1\rangle}$
\item for all ordinals $i\leq j$ less than
the length of $\eta$, we have $p_{\eta\restriction i}\leq p_\eta$,
and  there exist $f_{\eta\restriction
i,\eta}\in\Aut_{M_{\eta\restriction i}}\C$ such that $f_{\eta\restriction
i,\eta}(a_{\eta\restriction i})=a_\eta$ and $f_{\eta\restriction
i,\eta}=f_{\eta\restriction j,\eta}\circ f_{\eta\restriction
i,\eta\restriction j}$
\item $p_{\langle\rangle}=p$
\item $M_0=M$.
\end{enumerate}

Suppose that $M_i$ and $p_\eta\in\gaS(M_i)$ have been defined.  Since
$M_i$ is isomorphic to $M$ (by Fact \ref{unique limits}), our
assumption implies that $p_\eta$ cannot be quasi-minimal.  So we may fix
an extension
$N$ of
$M_i$ and two distinct big extensions of $p_\eta$ to $N$.  Let
$a'_{\eta\conc\langle 0\rangle}$ and
$a'_{\eta\conc\langle 1\rangle}$ realize these big extensions
and let
$M'_1\in\K_\mu$ be some extension of $N$ containing both
$a'_{\eta\conc\langle 0\rangle}$ and
$a'_{\eta\conc\langle 1\rangle}$.  Fix a
$(\mu,\omega)$-limit model over $N$ and call it $M_{i+1}$.
By the definition of big types, there are $a_{\eta\conc\langle
0\rangle}$ and $a_{\eta\conc\langle 1\rangle}$ realizing big extensions
of $\tp(a'_{\eta\conc\langle 0\rangle}/N)$ and
$\tp(a'_{\eta\conc\langle 1\rangle}/N)$, respectively.

For the limit stage of the construction notice that
$M_i:=\Union_{j<i}M_j$ is a limit model as guaranteed by condition
(\ref{get limit}).  For $\eta\in\sq i 2$ with $i$ a limit ordinal, we
choose $p_\eta$ to be some (there may be more than one) non-algebraic
extension of  the $p_{\eta\restriction j}$ for $j<i$. This is possible
by our construction of the $f_{\eta\restriction j,\eta}$'s and Lemma
\ref{Baldwin's coherence}.  This lemma also gives us
the required $f_{\eta\restriction j,\eta}$'s.

To see that this construction is enough, let $i$ be the first
ordinal $<\mu^+$ such that $2^i>\mu$.  Then, $\{p_\eta\in\gaS(M_i)\mid
\eta\in\sq i 2\}$ witnesses that $\K$ is not Galois-stable in $\mu$.

\end{proof}

We need an extension property for quasi-minimal types in order to find
deep-rooted quasi-minimal types.
\begin{proposition}[Extension Property for Quasi-Minimal
Types]\label{exist min ext}  Let $\K$ be  categorical in some
$\lambda>\chi$.  Let $\mu$ be such that
$\LS(\K)\leq\mu\leq\lambda$. If $p\in \gaS(M)$ is quasi-minimal and $M$ is
a
$(\mu,\sigma)$-limit model for some limit ordinal
satisfying
$\LS(\K)\leq\sigma<\mu^+$, then for every $M'\in\K_{\leq\lambda}$
extending $M$, there is a quasi-minimal $q\in \gaS(M')$ such that $q$
extends
$p$.  Furthermore if $p$ does not $\mu$ split over some $N$, then $q$ can
be chosen so that $q$ does not $\mu$-split over $N$.
\end{proposition}

\begin{proof}
Without loss of generality $M'$ is a limit model over $M$.
Let $p\in\gaS(M)$ be quasi-minimal.
Since $M$ is  $(\mu,\sigma)$-limit model,  using Fact
\ref{Baldwin non-split ext}, we can find a
proper submodel
$N\prec_{\K}M$ of cardinality $\mu$ such that for every
$M'\in\K_{\leq\lambda}$ there exists $q\in \gaS(M')$ extending $p$ such
that $q$ does not $\mu$-split over $N$.
Suppose for the sake of contradiction that $q$ is not quasi-minimal.  Then
tameness  and
Proposition \ref{monotonicity of
minimal types} tells us that $q$ must be algebraic.  Let $a\in M'$
realize $q$ and $M^a\in\K_\mu$ contain $a$ with
$M\prec_{\K}M^a\prec_{\K}M'$.
Then $q\restriction M^a$ is also algebraic.
However, since $q\restriction M^a$ does not
$\mu$-split over $N$ and extends $p$, by Corollary  \ref{non
split ext of non alg} we see that
$q\restriction M^a$ is not-algebraic.  This gives us a contradiction.
\end{proof}

\begin{remark}\label{ext for min}
Proposition \ref{exist min ext} holds for minimal types as well.  Simply
replace quasi-minimal with minimal in the proof.  This will be used in
the last section of the paper.
\end{remark}

Propositions \ref{exist min ext} and \ref{exist
minimal}  are key to get
the existence of deep-rooted quasi-minimal types.

\begin{proposition}[Existence of deep-rooted quasi-minimal
types]\label{exist strongly rooted} Let $\K$ be categorical in some
$\lambda>\chi$.
Then  for every
$M'\in\K_{\lambda}$, there exists a deep-rooted quasi-minimal
$q\in\gaS(M')$.  Furthermore, if $M\prec_{\K}M'$ is a limit model of
cardinality $\mu$ with $\chi\leq\mu<\lambda$ and $p\in\gaS(M)$ is
quasi-minimal, then we can find $q\in\gaS(M')$ a  deep-rooted
quasi-minimal extension of $p$ with root $M$.

\end{proposition}

\begin{proof}
Fix $\mu$ with $\chi\leq\mu<\lambda$.  Notice that by Fact \ref{cat
implies stable}, $\K$ is Galois-stable in $\mu$.
 Choose $M\in\K_\mu$ to be some $\K$-substructure of $M'$.  Since $\K$ is
stable in $\mu$ and categorical in
$\lambda$, we may take $M$ to be a $(\mu,\sigma)$-limit model for some
limit ordinal $\sigma$ with $\LS(\K)\leq\sigma<\mu^+$.
  By Proposition
\ref{exist minimal} and monotonicity of quasi-minimal types, we can choose
$M$ such  that
 there is
a quasi-minimal type $p\in \gaS(M)$.
Then by
Proposition \ref{exist min ext}, there exists a quasi-minimal
$q\in \gaS(M')$ extending $p$.  $q$ is rooted with root $M_i$.

\end{proof}

\section{Vaughtian Pairs}
We will show that for deep-rooted quasi-minimal types, there are no true
Vaughtian-pairs.  This is a variation of the result in \cite{GrVa1} that
for rooted minimal types there are no Vaughtian-pairs.

\begin{definition}
Let $M\in\K$ be a limit model and $p\in\gaS(M)$ non-algebraic. Fix
$\mu\geq\|M\|$.
\begin{enumerate}
\item A pair of models $(N_0,N_1)$ is said to be a
\emph{$(p,\mu)$-Vaughtian pair} provided that $N_0,N_1$ both have
cardinality $\mu$ and $M\preceq_{\K}N_0\precneqq_{\K}N_1$ with no
$c\in N_1\backslash N_0$ realizing $p$.
\item A $(p,\mu)$-Vaughtian pair $(N_0,N_1)$ is a \emph{true
$(p,\mu)$-Vaughtian pair} iff $N_0$ and $N_1$ are both limit
models.
\end{enumerate}
\end{definition}


The ubiquity of the assumption
of categoricity in a successor cardinal in the literature concerning
Conjecture \ref{shelah's cat conj} can be explained by the proof of the
following central result.  The result uses a classical
Vaughtian-construction in the spirit of Morley's work,
and it appears in
\cite{Sh 394} as Claim
$(*)_8$ of Theorem 9.7.

\begin{fact}\label{no lambda,lambda vps}
Assume that $\K$ is categorical in some $\lambda^+>\LS(\K)^+$.  Then for
every limit model $M\in\K_{\leq\lambda}$ and every minimal type $p\in
\gaS(M)$, there are no true
$(p,\lambda)$-Vaughtian pairs.
\end{fact}

Using the fact that all saturated models are limit models; that the union
of an increasing chain of saturated models is saturated (Claim 6.7 of
\cite{Sh 394}) and Fact
\ref{big is non-alg},  the same argument for Fact \ref{no lambda,lambda
vps} can be carried out to yield the following proposition.
\begin{proposition}\label{no lambda,lambda qm vps}
Assume that $\K$ is categorical in some $\lambda^+>\LS(\K)^+$.  Then for
every limit model $M\in\K_{\leq\lambda}$ and every quasi-minimal type
$p\in
\gaS(M)$, there are no true
$(p,\lambda)$-Vaughtian pairs.
\end{proposition}

Notice that the previous argument works only when
$\lambda$ is strictly larger than $\LS(\K)$.  We will come back to this
issue in Section \ref{s:lambda,lambda+} and deal with the special case
in which $\LS(\K)=\chi=\lambda$ and $\K$ is categorical in both $\lambda$
and
$\lambda^+$.

The following Vaughtian-pair transfer theorem is a relative of Theorem
3.3 of
\cite{GrVa1}:

\begin{theorem}\label{vp transfer}
Suppose that $\K$ is categorical in some $\lambda^+>\LS(\K)$.
Let $p$ be a deep-rooted quasi-minimal type over a model $M$ of
cardinality
$\lambda^+$. Fix a root $N\prec_{\K}M$ of cardinality $\lambda$, with
$p\restriction N$ quasi-minimal. If
$\K$ has a
$(p,\lambda^+)$-Vaughtian pair, then
there is a true $(p\restriction N,\lambda)$-Vaughtian pair
$(N_0,N_1)$.

\end{theorem}

\begin{proof}
Suppose that $(N^0,N^1)$ form a $(p,\lambda^+)$-Vaughtian pair.  By
categoricity, we know that $N^0$ and $N^1$ are both saturated.

   Let $C$ denote the set of all
realizations of
$p\restriction N$ inside $N^1$.  Fix
$a\in N^1\backslash N^0$.

We now construct $\langle
N^0_i,N^1_i\in\K_\lambda\mid i<\lambda^+\rangle$ satisfying the following:

\begin{enumerate}

\item $N^0_0=N$

\item $N^\ell_i\precneqq_{\K}N^\ell$ for $\ell=0,1$

\item  the sequences $\langle N^0_i\mid i<\lambda^+\rangle$ and
$\langle N^1_i\mid i<\lambda^+\rangle$ are both
$\prec_{\K}$-increasing and continuous

\item\label{get limit cond} $N^\ell_{i+1}$ is a limit model over
$N^\ell_i$ for
$\ell=0,1$

\item $a\in N^1_i\backslash N^0_i$ and

\item \label{put c in} $C_i:=C\bigcap N^1_i\subseteq N^0_{i+1}$.


\end{enumerate}

The construction follows from the fact that both $N^0$ and $N^1$ are
saturated and homogeneous and the following:

\begin{claim}\label{no new realizations in N1}
If $d\in N^1$ realizes $p\restriction N^0_0$,
then $d\in N^0$.  Thus
$C\subseteq N^0$.

\end{claim}

\begin{proof}[Proof of Claim \ref{no new realizations in N1}]
Suppose that $d\in N^1\backslash N^0$  realizes $p\restriction N^0_0$.
Because $N^0$ is saturated, $\tp(d/N^0)$ is not only non-algebraic, it is
a big extension of
$p\restriction N^0_0$.  Since $p\restriction N^0_0$ is quasi-minimal, we
have that
$\tp(d/M)=p$.  Since
$(N^0,N^1)$ form a $(p,\lambda^+)$-Vaughtian pair, it must be the
case that
$d\in N^0$, contradicting our choice of $d$.

\end{proof}

The construction is enough:
Define
 $$E:=\left\{\begin{array}{c|l}{\delta<\lambda^+}& \delta\text{ is a limit
ordinal},\\
&\text{for all }i<\delta\text{ and }x\in N^1_i,\\
&\text{if there exists }
j<\lambda^+\text{ such that }x\in C_j, \\ &\text{then there exists }
j<\delta, \text{ such that }x\in C_j\end{array}\right\}.$$  Notice that
$E$ is a club. (We only use the fact that $E$ is non-empty.)  Fix
$\delta\in E$.

\begin{claim}\label{force vp claim}
For every $c\in N^1_\delta\cap C$, we have $c\in N^0_\delta$.

\end{claim}

\begin{proof}[Proof of Claim \ref{force vp claim}]
Since $\langle
N^1_i \mid i<\lambda^+\rangle$ is continuous, there is $i<\delta$ such
that $c\in N^1_i$.  Thus by the definition of $E$, there is a
$j<\delta$ with $c\in C_j$.  By condition $(\ref{put c in})$ of the
construction, we would have put $c\in N^0_{j+1}\prec_{\K}N^0_\delta$.

\end{proof}

%




Notice that $N^1_\delta\neq N^0_\delta$ since $a\in N^1_\delta\backslash
N^0_\delta$.
Thus
Claim \ref{force vp claim}
allows us the conclude that we have constructed a $(p\restriction
N,\lambda)$-Vaughtian pair.
We complete the proof by observing that
 condition (\ref{get limit cond}) of
the construction and our choice of a limit ordinal $\delta$ imply  that
both
$N^0_\delta$ and
$N^1_\delta$ are limit models.

\end{proof}

\begin{remark}\label{close vp remark}
The same proof of Theorem \ref{vp transfer} works with minimal in place of
quasi-minimal.  Thus for $p$ a minimal type, notice that for $(N_0,N_1)$
in the conclusion of  Theorem \ref{vp transfer} we have that $(N_0,N_1)$
is a true
$(p\restriction N_0,\lambda)$-Vaughtian pair and $p\restriction N_0$ is
minimal.  Furthermore $N_0$ is a limit model over $N$.  This extra
information will be used in Section \ref{s:lambda,lambda+}.
\end{remark}

\begin{corollary}\label{no vp thm}
Let $\lambda>\LS(\K)$.
If
$\K$ is categorical in  $\lambda^+$ and $p$ is
a deep-rooted quasi-minimal type over a model of cardinality $\lambda^+$,
then there are no
$(p,\lambda^{+})$-Vaughtian pairs.
\end{corollary}
\begin{proof}
Suppose that $(N^0,N^1)$ is a $(p,\lambda^+)$-Vaughtian pair and that $N$
is both a root of $p$ (with $p\restriction N$ quasi-minimal) and a limit
model of size
$\lambda$.  Then by Theorem
\ref{vp transfer}, there is
a true $(p\restriction N,\lambda)$-Vaughtian pair.    This
contradicts Proposition
\ref{no lambda,lambda qm
vps}.

\end{proof}

\begin{corollary}\label{many realizations}
Let $\lambda>\LS(\K)$.
If
$\K$ is categorical in $\lambda^+$, then every deep-rooted quasi-minimal
type over a model $N$ of cardinality $\lambda^+$ is realized
$\lambda^{++}$ times in every model of cardinality $\lambda^{++}$
extending $N$.

\end{corollary}

\begin{proof}
Suppose $M\in\K_{\lambda^{++}}$ realizes $p$ only $\alpha<\lambda^+$
times.

Let $A:=\{a_i\mid i<\alpha\}$ be an enumeration of the realizations of $p$
in $M$.  We can find $N_0\in\K_{\lambda^+}$ such that
$N\Union A\subseteq N_0\prec_{\K}M$.  Since $M$ has cardinality
$\lambda^{++}$, we can
find $N_1\in\K_{\lambda^+}$  such that
$N_0\precneqq_{\K}N_1\prec_{\K}M$.  Then $(N_0,N_1)$ forms a $(p,
\lambda^+)$-Vaughtian pair contradicting
Corollary \ref{no vp thm}.

\end{proof}

\section{The Main Result}\label{s:main}

Now that we have established the existence of deep-rooted quasi-minimal
types with no Vaughtian pairs, we proceed as in \cite{GrVa1} to transfer
categoricity upwards using the following result which is a variation of
Theorem 4.1 of \cite{GrVa1}.

\begin{theorem}\label{many realizations of
minimal imply sat}
Let $\lambda\geq\chi$.  Suppose
$M_0\in\K_\lambda$ and
$r\in\gaS(M_0)$ is a quasi-minimal type  such that $\K$ has no
$(r,\lambda)$-Vaughtian pairs.

Let $\alpha$ be an ordinal $<\lambda^+$ such that
$\alpha=\lambda\cdot\alpha$.   Suppose $M\in\K_{\lambda}$ has a
resolution
$\langle M_i\in\K_{\lambda}\mid i<\alpha\rangle$  such that for
every
$i<\alpha$, there is
$c_i\in M_{i+1}\backslash M_i$ realizing $r$.  Then $M$ is saturated over
$M_0$.  Moreover if $\K$ is Galois-stable in $\lambda$, then $M$ is a
$(\lambda,\alpha)$-limit model over $M_0$.
\end{theorem}

\begin{proof}
At the referee's request we have included a proof of this result.
Let $r$, $M_0$, $M$ and $\langle M_i\mid i<\alpha\rangle$ be as in the
statement of the theorem.
Let $p\in\gaS(M_0)$ be given.  We will show that $M$ realizes $p$.

First, fix $M'$ an extension of $M_0$ of cardinality $\lambda$ realizing
$p$.  It is enough to construct an isomorphism between $M$ and some
extension of $M'$.  We build such an extension and isomorphism by
inductively defining increasing and continuous sequences $\langle M'_i\mid
i<\alpha\rangle$ and
$\langle h_i\mid i<\alpha\rangle$ so that $h_i:M_i\rightarrow M'_i$.
During this construction we also fix $\langle a_{\lambda\cdot i+j}\mid
j<\lambda\rangle$ an enumeration (possibly repeating) of the realizations
of $r$ inside $M_i$.  After stage
$\beta=\lambda\cdot i+j$ of the construction,
we  require that $a_{\beta}\in
h_{\beta+1}(M_{\beta +1})$.

To see that such a construction is possible, let us examine the successor
case.  The base and limit stages of the construction are routine to carry
out.  Suppose that we have defined $M'_i$ and $h_i$ and that we have
fixed an enumeration $\langle a_{\lambda\cdot k+j}\mid
j<\lambda\rangle$ of all realizations of $r$ in $M'_k$ for each $k\leq
i$.  By properties of ordinal arithmetic, there is exactly one pair $j,k$
with $k\leq i$ for which
$i+1=\lambda\cdot k+1$.
If $a_i$ is already in $h_{i}(M_i)$ there is nothing to do but extend
$h_i$ to include $M_{i+1}$ in its domain and choose an appropriate
$M'_{i+1}$ containing $h_{i+1}(M_{i+1})$.
When $a_i\notin h_i(M_i)$, more care is needed.
The important thing to notice here is
that in this case,
$\tp(a_i/h_i(M_i))$ is a non-algebraic extension of $r$.  By the
quasi-minimality of $r$, we know that regardless of which extension
$\check h$ of
$h_i$ to an automorphism of $\C$ that one would consider, we
have $\tp(\check h(c_i)/h_i(M_i))=\tp(a_i/h_i(M_i))$.  Thus we can choose
$\check h$ to be an automorphism of $\C$ extending $h_i$ so that
$\check h(c_i)=a_i$.  Now define $h_{i+1}:=\check h\restriction M_{i+1}$
and choose an appropriate extension $M'_{i+1}$ of $M'_i$ containing the
image of
$M_{i+1}$  under $h_{i+1}$.

Once we have completed the construction outlined above, the issue of
whether or not $h:=\Union_{i<\alpha}h_i$ is an isomorphism between $M$
and $\Union_{i<\alpha}M'_i$ remains to be addressed.  First notice that by
our assumption that
$\alpha=\lambda\cdot\alpha$, if $a\in M'$ realizes $r$, then at some
stage in the construction, we would have put $a$ into the range of $h$.
Therefore, if $h$ were not an isomorphism,  $h(M)$ and $M'$ would form a
$(r,\lambda)$-Vaughtian pair contradicting our hypothesis on $r$.

If in addition to the hypothesis given, we assume that $\K$ is
Galois-stable in $\lambda$, we could conclude that $M$ is a
$(\lambda,\alpha)$-limit model by altering the construction.
At stage $i$ of the construction we choose $M'_{i+1}$ as above,
only now require that
$M'_{i+1}$ to be universal over $M'_i$.
\end{proof}

Using Theorem \ref{many realizations of
minimal imply sat}, we are able to transfer categoricity from $\lambda$
to $\lambda^+$ by showing that every model of cardinality $\lambda^+$ is
saturated:

\begin{theorem}
Suppose that $\K$ has arbitrarily large models, is $\chi$-tame and
satisfies the amalgamation and joint embedding properties. Let $\lambda$
be such that
$\lambda>\LS(\K)$ and $\lambda\geq\chi$. If
$\K$ is categorical in
$\lambda^+$ then
$\K$ is categorical in all $\mu\geq\lambda^{+}$.
\end{theorem}

\begin{proof}

First we prove that $\K$ is categorical in $\lambda^{++}$ by establishing
that every model $N$ of cardinality $\lambda^{++}$ is saturated.
Let $M\prec_{\K}N$ have cardinality $\lambda^+$.  We will show that $N$
realizes every type over $M$.  First notice that Proposition \ref{exist
strongly rooted} and categoricity in
$\lambda^+$ guarantees that there exists  a deep-rooted quasi-minimal
$r\in\gaS(M)$. By Corollary \ref{many
realizations}, we know that $N$ realizes $r$ $\lambda^{++}$-times.

 Let $\alpha<\lambda^+$ be such that
$\alpha=\lambda^+\cdot\alpha$. By the
Downward-L\"{o}wenheim Skolem Axiom of AECs, we can construct a
$\prec_{\K}$-increasing and continuous chain of models
$\langle M_i\prec_{\K}N\mid i<\alpha\rangle $ such that $M=M_0$
for every $i<\alpha$, we can fix $a_i\in M_{i+1}\backslash M_i$ realizing
$r$.  This construction is possible since there are $\lambda^{++}$-many
realizations of $r$ from which to choose.  By Fact \ref{many realizations of
minimal imply sat},
$\Union_{i<\alpha}M_i$ realizes every type over $M$.

We have explained that categoricity in $\lambda^+$
implies categoricity in $\lambda^+$ and $\lambda^{++}$.  Now, an
application of Fact \ref{going up grva} provides categoricity in all
larger cardinalities.

\end{proof}

A combination of our upward result and Shelah's downward result from
\cite{Sh 394} yields
\begin{theorem}\label{abstract thm}
Let $\K$ be a $\chi$-tame abstract elementary class  satisfying the
amalgamation and joint embedding properties.  If $\K$ is categorical in
 $\lambda^+$ for some $\lambda >\max\{\LS(\K),\chi\}$, then $\K$ is
categorical in all
$\mu\geq\min\{\lambda^+,\beth_{(2^{\Hanf(\K)})^+}\}$.
\end{theorem}
  It remains open whether
or not categoricity in $\lambda^+$ implies categoricity in $\lambda^{++}$
for the special case where
$\aleph_0<\LS(\K)=\chi=\lambda$.  For this case, a substitute for
Fact
\ref{no lambda,lambda vps} is missing.  We will provide some partial
results  concerning this problem in the following section.

\section{Categoricity in $\LS(\K)$ and
$\LS(\K)^+$}\label{s:lambda,lambda+}

In this section, we examine an abstract
elementary class which is categorical in both $\lambda$ and
$\lambda^+$ and $\lambda=\LS(\K)=\chi$.  We assume the class has no
maximal models and satisfies the amalgamation and joint embedding
properties.  This is motivated by questions of John Baldwin and Olivier
Lessmann concerning perceived limitations of
\cite{GrVa1}.  From these assumptions, we derive categoricity in all
$\mu\geq\LS(\K)$.   The difficulty in working with a class that is
categorical in
$\LS(\K)^+$ is that there are no saturated models of cardinality
$\LS(K)$.  However, from stability we do have limit models of cardinality
$\LS(\K)$, and in this section we have an extra categoricity assumption
which tells us that all models of cardinality $\LS(\K)$, while not
saturated, are limit models.  This allows us to use minimal types instead
of quasi-minimal types.

We begin with a replacement for Fact \ref{no lambda,lambda vps}.
\begin{theorem}\label{no vps for lambda,lambda+}
Assume that $\K$ is categorical in $\lambda$ and $\lambda^+$ with
$\lambda=\LS(\K)=\chi$.  Then for every limit model $M\in\K_{\lambda}$
there is a minimal type
$p\in
\gaS(M)$, such that there are no true
$(p,\lambda)$-Vaughtian pairs of the form $(N_0,N_1)$ with $M=N_0$.
\end{theorem}

\begin{proof}
Suppose every minimal type over a limit model had a true Vaughtian pair.
Let $M$ be a limit model of cardinality $\mu$ and fix $p\in\gaS(M)$
minimal with true Vaughtian pair $(M,N_1)$ where
$N_1\in\K_\lambda$.  We can construct a $\prec_{\K}$-increasing and
continuous chain $\langle N_i\mid i<\lambda^+\rangle $ of limit models
such that for each $i<\lambda^+$
\begin{enumerate}
\item $N_0=M$
\item $N_i\in\K_\lambda$
\item $N_i$ is a limit model and
\item\label{no new real cond} no $a\in N_{i+1}\backslash N_i$ realizes
$p$.
\end{enumerate}
Suppose
$i$ is a limit ordinal and that we have defined $N_j$ for all $j<i$.
Let $N_i:=\Union_{j<i}N_j$.  By categoricity in $\lambda$ we know that
$N_i$ must be a limit model (but it may not be a limit model over $M$).

For the successor step of the construction,  suppose that $N_i$ has been
defined. Since
$M$ is a limit model, we can find
$p_i\in\gaS(N_i)$ a unique non-algebraic extension of $p$ (by Remark
\ref{ext for min}).  Since $N_i$ is a limit model and $p_i$ is a
minimal type, by our assumption it must be the case that there is
$N_{i+1}$ a limit model extending $N_i$ which together with $N_{i}$
forms a
$(p_i,\lambda)$-Vaughtian pair.  Since no $a\in N_{i+1}\backslash N_i$
realizes $p_i$, we can conclude by the minimality of $p$ that condition
(\ref{no new real cond}) holds.

To see why the construction is enough to get a contradiction, let
$N_{\lambda^+}:=\Union_{i<\lambda^+}N_i$.  From
condition (\ref{no new real cond}) of the construction, we
find that $N_{\lambda^+}$ does not realize
$p$.  Thus $N_{\lambda^+}$ is not saturated, which contradicts
categoricity in
$\lambda^+$.
\end{proof}

We now prove a slight variation of Corollary \ref{no vp thm}.
\begin{corollary}\label{no vp cor}
Let $\lambda$ be as in Theorem \ref{no vps for lambda,lambda+}.
For every $M\in\K_{\lambda^+}$, there is $q\in\gaS(M)$,
a deep-rooted minimal type with no
$(q,\lambda^{+})$-Vaughtian pairs.
\end{corollary}

\begin{proof}
Let $M\in\K_{\lambda^+}$ be given.  Fix $N\prec_{\K}M$ a limit model of
cardinality $\lambda$.  By Theorem \ref{no vps for lambda,lambda+}, we
can choose a minimal $p\in\gaS(N)$ such that there are no true
$(p,\lambda)$-Vaughtian pairs.  By Proposition \ref{exist strongly
rooted}, we can extend $p$ to a deep-rooted minimal type $q\in\gaS(M)$.
Suppose that $N^0,N^1$ form a $(q,\lambda^+)$-Vaughtian pair.
Then
Theorem
\ref{vp transfer} and Remark \ref{close vp remark} tell us that there
are limit models $N_0,N_1$ with $N\prec_{\K}N_0\prec_{\K}N^0$ with
$(N_0,N_1)$ a $(p,\lambda)$-Vaughtian pair and $N_0$ a limit model over
$N$.  Furthermore, we have that $(N_0,N_1)$ form a $(q\restriction
N_0,\lambda)$-Vaughtian pair.

We will now show that by our choice of $p$
such $(q\restriction N_0,\lambda)$-Vaughtian pairs cannot exist.
  Since
$N$ is a limit model, we can find a resolution
$\langle N^+_i\mid i<\omega\rangle$ of $N$ such that $N^+_{i+1}$ is
universal over
$N^+_i$.  By Fact \ref{Baldwin
non-split ext}, there is $i<\omega$ such that $p$ does not
$\lambda$-split over $N^+_i$.  Observe that $N$ is a limit model over
$N^+_{i}$. Additionally, since $N_0$ is a limit model over $N$ it is
also a limit model over
$N^+_{i}$.  Then, $N$ and $N_0$ are isomorphic over $N^+_i$.  Let
$f:N\cong N_0$ with $f\restriction N^+_i=\id_{N^+_i}$.  Since there are no
$(p,\lambda)$-Vaughtian pairs with $N$ as the first model in the pair,
there are no $(f(p),\lambda)$-Vaughtian pairs with $N_0$ as the first
model in the pair.  By invariance and our choice of $N^+_i$, we have that
$f(p)$ does not $\mu$-spit over $N^+_i$.  This implies that
$f(p)\geq p$, otherwise $f^{-1}$ would witness that $f(p)$
$\lambda$-splits over $N^+_i$.  Now we have that $f(p)$ and
$q\restriction N_0$ are both non-algebraic extensions of
$p$ to $N_0$.  By minimality of $p$, $f(p)=q\restriction N_0$ and we can
conclude that there are no $(q\restriction N_0)$-Vaughtian pairs with
$N_0$ as the first model of the pair.  This gives us a contradiction and
completes the proof.
\end{proof}

Corollary \ref{no vp cor} is enough to carry out the argument of
Corollary \ref{many realizations} and the remaining arguments in Section
\ref{s:main}.  This allows us to conclude the second theorem in the
abstract, restated here:

\begin{theorem}
Let $\K$ be a $\LS(\K)$-tame abstract elementary class satisfying the
amalgamation and joint embedding properties with arbitrarily large
models.  If
$\K$ is categorical in both $\LS(K)$ and $\LS(\K)^+$, then $\K$ is
categorical in all $\mu\geq\LS(\K)$.
\end{theorem}

\end{document}